\documentclass{article}
\usepackage{amssymb}
\usepackage{amsfonts}

\usepackage{bm}

\usepackage[toc,page]{appendix}

\usepackage{cite}

\usepackage{etoolbox}
\usepackage{chngcntr}
\usepackage{amsmath}
\numberwithin{equation}{section}

\evensidemargin=\oddsidemargin
\newdimen\dummy
\dummy=\oddsidemargin
\title{The Legendre-Hadamard condition in Cosserat elasticity theory}
\author{Milad Shirani$^{1a}$, David J. Steigmann$^{1b}$ and Patrizio Neff$^{2}$}
\date{\empty}

\begin{document}

\maketitle

\begin{center}
%
%
%

$^{1}$Department of Mechanical Engineering

University of California

Berkeley, CA. 94720, USA

$^{a}$milad\_shirani@berkeley.edu, $^{b}$dsteigmann@berkeley.edu

\qquad \qquad 

$^{2}$Chair of Nonlinear Analysis and Modelling

Fakult\"{a}t f\"{u}r Mathematik

Universit\"{a}t Duisburg-Essen

45141 Essen, Germany

patrizio.neff@uni-due.de

\quad

20 May, 2020
\end{center}

\qquad \qquad 

\textit{Summary}: The Legendre-Hadamard necessary condition for energy
minimizers is derived in the framework of Cosserat elasticity theory.

\qquad

\section{Introduction}


Cosserat elasticity \cite{cosserat1909theorie} is enjoying a resurgence as a framework for
the modeling and analysis of scale effects in solids associated with the
presence of microstructure. Definitive modern treatments of the subject may
be found in \cite{toupin1964theories, reissner1975note, reissner1987further, truesdell2004non, neff2006existence, pietraszkiewicz2009natural, neff2015existence, lankeit2017integrability, neff2017real, neff2004geometrically}. Here we supplement this literature
with the relevant version of the Legendre-Hadamard necessary condition for
energy minimizers. Thus we effectively extend the Legendre-Hadamard
inequality of conventional elasticity theory \cite{truesdell2004non} to the Cosserat
theory.

We work in the setting of classical nonlinear Cosserat theory, according to
which the material comprising the considered body is endowed with
independent deformation and rotation fields, the former describing the
transplacements of material points as in conventional elasticity theory and
the latter the change in microstructural orientation as the configurations
of the body evolve.

Section 2 is devoted to a brief resum\'{e} of the basic theory for
equilibria. We regard the latter as those states that satisfy an appropriate
virtual work statement. Conditions under which this may be interpreted as a
stationarity condition for a potential energy functional are identified in
Section 3, and expressions for the first and second variations of this
energy are obtained. In Section 4 we present a detailed derivation, modelled
after that given in \cite{giaquinta2004weierstrass}, of the relevant Legendre-Hadamard
inequality. This proceeds from the notion that the second variation is
necessarily non-negative if an equilibrium state furnishes a minimum of the
potential energy. We conclude in Section 5 with an application of the
inequality to the particular strain-energy function proposed in \cite{neff2015existence}.

Concerning notation, bold face is used for vectors and tensors and a dot
interposed between bold symbols is used to denote the standard Euclidean
inner product. For example, if $\bm{A}$ and $\bm{B}$ are
second-order tensors, then their inner product is $\bm{A}\cdot\bm{B}=tr(\bm{AB}^{t}),$ where $tr(\cdot )$ is the trace
and the superscript $^{t}$ is used to denote the transpose. The induced norm
is $\left\Vert \bm{A}\right\Vert =\sqrt{\bm{A}\cdot \bm{A}}.$ We make
frequent use of the fact that $\bm{A}\cdot \bm{BC}=\bm{AC}^{t}\cdot \bm{B}.$
The symbol $\otimes $ identifies the standard tensor product of vectors. We
use $sym\bm{A}$ and $skew\bm{A}$ respectively to denote the
symmetric and skew parts of a tensor $\bm{A}$, and $dev(sym\bm{A)}$
to denote the deviatoric part of $sym\bm{A}.$ The axial vector of a skew
tensor $\bm{W}$ is denoted by $axl\bm{W}$ and defined by $axl\bm{%
W}\wedge \bm{v}=\bm{Wv}$ for any vector $\bm{v}.$ The symbols $\nabla $ and $Div$
respectively stand for the three-dimensional referential gradient and
divergence operators. For a fourth-order tensor $\mathcal{A}$, the notation $%
\mathcal{A}[\bm{B}]$ stands for the second-order tensor resulting from
the linear action of $\mathcal{A}$ on $\bm{B}$ (see \cite{truesdell2004non}, eq.
(7.10)). Its transpose $\mathcal{A}^{t}$ is defined by $\bm{B}\cdot 
\mathcal{A}^{t}[\bm{A}]=\bm{A}\cdot \mathcal{A}[\bm{B}],$ and $%
\mathcal{A}$ is said to possess major symmetry if $\mathcal{A}^{t}=\mathcal{A%
}$. The notation $G_{\bm{S}}$ stands for the second-order-tensor-valued
derivative of the scalar-valued function $G(\bm{S})$ with respect to the
second-order tensor variable $\bm{S}$. The second derivative is the
fourth-order tensor $G_{\bm{SS}}$; this possesses major symmetry if $G$
is twice differentiable. The second derivatives $G_{\bm{ST}}$ and $G_{%
\bm{TS}}$ of a twice differentiable scalar-valued function $G(\bm{%
S},\bm{T})$ satisfy $\bm{A}\cdot G_{\bm{ST}}[\bm{B}]=\bm{B}\cdot G_{%
\bm{TS}}[\bm{A}];$ accordingly, $G_{\bm{TS}}=(G_{\bm{ST}%
})^{t}.$ Finally, we use superposed dots to denote variational derivatives.
These are ordinary derivatives of one-parameter families of the varied
functions with respect to the parameter, evaluated at parameter value zero,
say, which we identify with an equilibrium state.


\section{Cosserat elasticity}


We present a brief outline of the equilibrium theory for the sake of
completeness.


\subsection{Kinematics and elasticity}


The relevant kinematical variables of a Cosserat continuum are a deformation
field $\bm{\chi} (\bm{X})$ and a rotation field $\bm{R}(\bm{X})$, where $%
\bm{X}$ is the position of a material point in a reference configuration 
$\kappa ,$ say. Of course these may depend on time, but such dependence is
not important for our purposes and is thus not made explicit. The
deformation and rotation fields are regarded as being independent in the
spirit of the conventional Cosserat theory (see \cite{truesdell2004non}, \S\ 98).

To model elasticity, we introduce an energy density $U(\bm{F},\bm{R},\nabla 
\bm{R};\bm{X}),$ per unit volume of $\kappa ,$ where $\bm{F}=\nabla 
\bm{\chi }$ is the deformation gradient and $\nabla \bm{R}$ is the
rotation gradient. In Cartesian index notation, these are 
\begin{equation}
\bm{F}=F_{iA}\bm{e}_{i}\otimes \bm{E}_{A},\quad\bm{R}%
=R_{iA}\bm{e}_{i}\otimes \bm{E}_{A}\quad \textrm{and}\quad \nabla 
\bm{R}=R_{iA,B}\bm{e}_{i}\otimes \bm{E}_{A}\otimes \bm{E}_{B}
\label{eq:2.1}
\end{equation}%
with%
\begin{equation}
F_{iA}=\chi _{i,A},  \label{eq:2.2}
\end{equation}%
where $(\cdot ),_{A}=\partial (\cdot )/\partial X_{A}$. Here $\{\bm{e}%
_{i}\}$ and $\{\bm{E}_{A}\}$ are fixed orthonormal bases associated with
Cartesian coordinates $x_{i}$ and $X_{A},$ where $x_{i}=\chi _{i}(X_{A}).$

We assume the strain energy to be Galilean-invariant and thus impose%
\begin{equation}
U(\bm{F},\bm{R},\nabla\bm{R};\bm{X})=U(\bm{QF},\bm{QR},\bm{Q}\nabla \bm{R};\bm{X}),  \label{eq:2.3}
\end{equation}%
where $\bm{Q}$ is an arbitrary spatially uniform rotation and $(\bm{Q%
}\nabla \bm{R})_{iAB}=(Q_{ij}R_{jA}),_{B}=Q_{ij}R_{jA,B}.$ The
restriction%
\begin{equation}
U(\bm{F},\bm{R},\nabla \bm{R};\bm{X})=W(\bm{E},\bm{\Gamma} ;\bm{X}), 
\label{eq:2.4}
\end{equation}%
where \cite{pietraszkiewicz2009natural}%
\begin{equation}
\bm{E=R}^{t}\bm{F}=E_{AB}\bm{E}_{A}\otimes \bm{E}_{B};\quad
E_{AB}=R_{iA}F_{iB},  \label{eq:2.5}
\end{equation}%
and%
\begin{equation}
\bm{\Gamma }=\Gamma _{DC}\bm{E}_{D}\otimes \bm{E}_{C};\quad
\Gamma _{DC}=\frac{1}{2}e_{BAD}R_{iA}R_{iB,C},  \label{eq:2.6}
\end{equation}%
with $W$ the reduced energy and $e_{ABC}$ the permutation symbol ($%
e_{123}=1, $ etc.), furnishes the necessary and sufficient condition for
Galilean invariance. Sufficiency is obvious, whereas necessity follows by
choosing $\bm{Q}=\bm{R}_{\mid x}^{t},$ where $x$ is the material point in
question, and making use of the fact that for each fixed $C\in \{1,2,3\},$
the matrix $R_{iA}R_{iB,C}$ is skew. This follows by differentiating $%
R_{iA}R_{iB}=\delta _{AB}$ (the Kronecker delta). The associated axial
vectors $\bm{\gamma }_{C}$ have components 
\begin{equation}
\gamma _{D(C)}=\frac{1}{2}e_{BAD}R_{iA}R_{iB,C},  \label{eq:2.7}
\end{equation}%
yielding \cite{pietraszkiewicz2009natural} 
\begin{equation}
\bm{\Gamma }=\bm{\gamma }_{C}\otimes \bm{E}_{C},  \label{eq:2.8}
\end{equation}%
and so $\bm{\Gamma }$ - the \textit{wryness tensor} - is isomorphic to
the Cosserat strain measure $\bm{R}^{t}\nabla \bm{R}.$ The strain
measures $\bm{E}$ and $\bm{\Gamma }$ are generally non-symmetric.

We note that the considerations of  \cite{neff2015existence} and \cite{lankeit2017integrability} are based
on strain measures that differ from those adopted here. However, in these
works it is demonstrated that the various sets of measures adopted therein
are equivalent to those used in the present work.

Henceforth we assume $W$ to be a continuous function of $\bm{X}$ and
twice continuously differentiable with respect to $\bm{E}$ and $\bm{%
\Gamma}.$


\subsection{Virtual power and equilibrium}


We define equilibria to be states that satisfy the virtual-power statement%
\begin{equation}
\dot{S}=P,  \label{eq:2.9}
\end{equation}%
where $P$ is the virtual power of the loads acting on the body, the explicit
form of which is deduced below, 
\begin{equation}
S=\int_{\kappa }Udv  \label{eq:2.10}
\end{equation}%
is the total strain energy, and, here and henceforth, superposed dots
identify variational derivatives. Thus, by the chain rule, 
\begin{equation}
\dot{U}=\dot{W}=\bm{\sigma} \cdot \bm{\dot{E}}+\bm{\mu} \cdot \bm{\dot{\Gamma}},%
  \label{eq:2.11}
\end{equation}%
where 
\begin{equation}
\bm{\sigma }=W_{\bm{E}}\quad \textrm{and}\quad \bm{\mu }=W_{%
\bm{\Gamma }}  \label{eq:2.12}
\end{equation}%
are evaluated at equilibrium, i.e., at states satisfying (\ref{eq:2.9}).

It follows from (\ref{eq:2.5}) that%
\begin{equation}
\bm{\dot{E}}=\bm{R}^{t}(\nabla \bm{u}-\bm{\Omega F}),\textrm{%
\quad where\quad }\bm{u}=\bm{\dot{\chi}}\textrm{\quad and\quad }\bm{\Omega} =%
\bm{\dot{R}R}^{t},  \label{eq:2.13}
\end{equation}%
where $\bm{\Omega }$ is an arbitrary skew tensor (see the Appendix).

Then,%
\begin{equation}
\bm{\sigma} \cdot \bm{\dot{E}}=\bm{R\sigma }\cdot \nabla \bm{u}-\bm{\Omega 
}\cdot skew(\bm{R\sigma F}^{t}).  \label{eq:2.14}
\end{equation}%
Let $\bm{\omega }=axl\bm{\Omega }$. If $\bm{\alpha }$ is a skew
tensor and $\bm{a}=axl\bm{\alpha ,}$ then it is easy to show that $%
\bm{\Omega} \cdot \bm{\alpha }=2\bm{\omega} \cdot \bm{a}.$ Further, $\bm{%
R\sigma F}^{t}=\bm{R\sigma E}^{t}\bm{R}^{t}$ and $skew(\bm{%
R\sigma E}^{t}\bm{R}^{t})=\bm{R}skew(\bm{\sigma E}^{t})\bm{R}%
^{t},$ yielding%
\begin{equation}
\bm{\sigma} \cdot \bm{\dot{E}}=\bm{R\sigma} \cdot \nabla \bm{u}-2axl[\bm{R}%
skew(\bm{\sigma E}^{t})\bm{R}^{t}]\cdot \bm{\omega }.  \label{eq:2.15}
\end{equation}

The reduction%
\begin{equation}
\bm{\dot{\Gamma}}=\bm{R}^{t}\nabla \bm{\omega }  \label{eq:2.16}
\end{equation}%
is somewhat more involved. Reference may be made to \cite{eugster2019continuum} for a
detailed derivation.

Accordingly, 
\begin{equation}
\bm{\mu} \cdot \bm{\dot{\Gamma}}=\bm{R\mu }\cdot \nabla \bm{\omega } 
\label{eq:2.17}
\end{equation}%
and on substituting (\ref{eq:2.11}), (\ref{eq:2.15})\ and (\ref{eq:2.16}) into (\ref{eq:2.9}) we obtain 
\begin{eqnarray}
P &=&\int_{\partial \kappa }[(\bm{R\sigma })\bm{\nu} \cdot \bm{u}+(\bm{R\mu})\bm{\nu} \cdot
\bm{\omega}]da  \nonumber \\
&&-\int_{\kappa }\{\bm{u}\cdot Div(\bm{R\sigma})+\bm{\omega }%
\cdot \lbrack Div(\bm{R\mu})+2axl(\bm{R}skew(\bm{\sigma E}^{t})%
\bm{R}^{t})]\}dv,  \label{eq:2.18}
\end{eqnarray}%
where $\bm{\nu }$ is the exterior unit normal to the (piecewise smooth)
surface $\partial \kappa $. The virtual power is thus of the form 
\begin{equation}
P=\int_{\partial \kappa }(\bm{t}\cdot \bm{u}+\bm{c}\cdot \bm{\omega})da+\int_{\kappa }(%
\bm{g}\cdot \bm{u}+\bm{\pi} \cdot \bm{\omega})dv,  \label{eq:2.19}
\end{equation}%
where $\bm{t}$ and $\bm{c}$ are densities of force and couple acting
on $\partial \kappa ,$ and $\bm{g}$ and $\bm{\pi }$ are densities of
force and couple acting in $\kappa .$

If there are no kinematical constraints; that is, if $\bm{u}$ and $%
\bm{\omega }$ can be chosen independently and arbitrarily, then, by the
Fundamental Lemma,%
\begin{equation}
\bm{g}=-Div(\bm{R\sigma})\quad \textrm{and}\quad \bm{\pi }=-Div(%
\bm{R\mu})-2axl[\bm{R}skew(\bm{\sigma E}^{t})\bm{R}%
^{t}]\quad \textrm{in}\quad \kappa ,  \label{eq:2.20}
\end{equation}%
whereas 
\begin{equation}
\bm{t}=(\bm{R\sigma})\bm{\nu}\textrm{\quad on\quad }\partial \kappa _{t}\quad 
\textrm{and\quad }\bm{c}=(\bm{R\mu})\bm{\nu}\textrm{\quad on\quad }\partial \kappa
_{c},  \label{eq:2.21}
\end{equation}%
where $\partial \kappa _{t}$ is a part of $\partial \kappa $ where position
is not assigned and $\partial \kappa _{c}$ is a part where rotation is not
assigned. We assume position to be assigned on $\partial \kappa \setminus
\partial \kappa _{t}$, so that $\bm{u}=\bm{0}$ there, and rotation to be
assigned on $\partial \kappa \setminus \partial \kappa _{c}$, where $\bm{%
\omega} =\bm{0}$. These, in addition to the degree of smoothness implied by the
foregoing reduction, are the admissibility conditions on $\bm{u}$ and $%
\bm{\omega .}$

Equations (\ref{eq:2.20}) and (\ref{eq:2.21}) are the equilibrium conditions for an elastic
Cosserat continuum.


\section{Conservative problems and potential energy}


We are concerned in this work with conservative problems for which a
potential energy is available. These are such that there exists a \textit{%
load potential }$L,$ say, whose variational derivative is identical to the
virtual power. Thus,%
\begin{equation}
\dot{L}=P  \label{eq:3.1}
\end{equation}%
and the potential energy is%
\begin{equation}
E=S-L,  \label{eq:3.2}
\end{equation}%
apart from an unimportant constant. Equilibria are thus seen to be those
states that render the potential energy stationary, i.e.,%
\begin{equation}
\dot{E}=0,  \label{eq:3.3}
\end{equation}%
for all admissible $\bm{u}$ and $\bm{\omega .}$


\subsection{Dead-load problems}


For the sake of simplicity and definiteness we confine attention to
dead-load problems with vanishing volumetric densities of force $\bm{g}$
and couple $\bm{\pi }$. These are characterized by load potentials of
the form%
\begin{equation}
L=\int_{\partial \kappa _{t}}\bm{t}\cdot \bm{\chi} da +\int_{\partial \kappa
_{c}}\bm{M} \cdot \bm{R}da  \label{eq:3.4}
\end{equation}%
in which $\bm{t}$ and $\bm{M}$ respectively are assigned
configuration-independent vector and tensor fields. Here $\bm{t}$ is as
in (\ref{eq:2.21})$_{1}$, and the (configuration dependent) couple traction in (\ref{eq:2.21})$%
_{2}$ is%
\begin{equation}
\bm{c}=2axl[skew(\bm{MR}^{t})].  \label{eq:3.5}
\end{equation}

The first variation of the energy is%
\begin{equation}
\dot{E}=\int_{\kappa }(\bm{R}W_{\bm{\Gamma }}\cdot \nabla \bm{%
\omega }+\bm{R}W_{\bm{E}}\cdot \nabla \bm{u}-\bm{R}W\bm{%
_{\bm{E}}F}^{t}\cdot \bm{\Omega })dv-\int_{\partial \kappa _{t}}%
\bm{t}\cdot \bm{u}da-\int_{\partial \kappa _{c}}\bm{c}\cdot \bm{\omega} da, 
\label{eq:3.6}
\end{equation}%
and vanishes if and only if the state $\{\bm{\chi} ,\bm{R}\}$ is equilibrated.


\subsection{The second variation at equilibrium}


To secure an expression for the second variation, we define $\bm{v}=\bm{\ddot{%
\chi}},$ and note, from (\ref{eq:2.13})$_{3},$ that 
\begin{equation}
\bm{\ddot{R}}=\bm{\Phi R}+\bm{\Omega}^{2}\bm{R},  \label{eq:3.7}
\end{equation}%
where $\bm{\Omega }$ is defined in (\ref{eq:2.13})$_{3}$ and $\bm{\Phi }$ is
an arbitrary skew tensor (see the Appendix).

On taking a further variation of (\ref{eq:3.6}), after some effort we obtain%
\begin{eqnarray}
\ddot{E} &=&\int_{\kappa }(\bm{R}W_{\bm{\Gamma }}\cdot \nabla 
\bm{\varphi}+\bm{R}W_{\bm{E}}\cdot \nabla \bm{v}-\bm{R}W%
\bm{_{\bm{E}}F}^{t}\cdot \bm{\Phi })dv-\int_{\partial \kappa
_{t}}\bm{t}\cdot \bm{v}da-\int_{\partial \kappa _{c}}\bm{c}\cdot \bm{\varphi}%
da  \nonumber \\
&&+\int_{\kappa }[\bm{\Omega R}W_{\bm{\Gamma }}\cdot \nabla \bm{%
\omega }+\bm{\Omega R}W_{\bm{E}}\cdot \nabla \bm{u}-\bm{%
\Omega R}W\bm{_{\bm{E}}F}^{t}\cdot \bm{\Omega }-\bm{R}W%
\bm{_{\bm{E}}}(\nabla \bm{u)}^{t}\cdot \bm{\Omega ]}dv 
\nonumber \\
&&+\int_{\kappa }[\bm{R}(W_{\bm{E}})^{\cdot }\cdot \nabla \bm{u}+%
\bm{R}(W_{\bm{\Gamma }})^{\cdot }\cdot \nabla \bm{\omega -R}(W_{%
\bm{E}})^{\cdot }\bm{F}^{t}\cdot \bm{\Omega ]}dv-\int_{\kappa
_{c}}\bm{MR}^{t}\cdot \bm{\Omega }^{2}da,  \label{eq:3.8}
\end{eqnarray}%
where $\bm{\varphi }=axl\bm{\Phi ,}$ and, by the chain rule,%
\begin{equation}
(W_{\bm{E}})^{\cdot }=W_{\bm{EE}}[\bm{\dot{E}}]+W_{\bm{%
E\Gamma }}[\bm{\dot{\Gamma}}]\textrm{\quad and\quad }(W_{\bm{\Gamma }%
})^{\cdot }=W_{\bm{\Gamma E}}[\bm{\dot{E}}]+W_{\bm{\Gamma \Gamma 
}}[\bm{\dot{\Gamma}}]  \label{eq:3.9}
\end{equation}%
with $\bm{\dot{E}}$ and $\bm{\dot{\Gamma}}$ given by (\ref{eq:2.13})$_{1}$
and (\ref{eq:2.16}), respectively. Here $\bm{v}$ vanishes on $\partial \kappa
\diagdown \partial \kappa _{t}$ and $\bm{\varphi }$ vanishes on $%
\partial \kappa \diagdown \partial \kappa _{c}.$

If the state $\{\bm{\chi} ,\bm{R}\}$ is equilibrated then the first line of
(\ref{eq:3.8})\ vanishes by (\ref{eq:3.3}) and (\ref{eq:3.6}). The second variation at equilibrium
becomes%
\begin{eqnarray}
\ddot{E} &=&\int_{\kappa }\{\bm{R}^{t}\nabla \bm{u}\cdot W_{\bm{%
EE}}[\bm{R}^{t}\nabla \bm{u}]+\bm{R}^{t}\nabla \bm{u}\cdot
W_{\bm{E\Gamma }}[\bm{R}^{t}\nabla \bm{\omega }]+\bm{R}%
^{t}\nabla \bm{\omega }\cdot W_{\bm{\Gamma E}}[\bm{R}^{t}\nabla 
\bm{u}]+\bm{R}^{t}\nabla \bm{\omega }\cdot W_{\bm{\Gamma
\Gamma }}[\bm{R}^{t}\nabla \bm{\omega }]\}dv  \nonumber \\
&&+\int_{\kappa }F(\nabla \bm{u},\nabla \bm{\omega} ,\bm{\Omega})%
dv-\int_{\kappa _{c}}\bm{MR}^{t}\cdot \bm{\Omega }^{2}da, 
\label{eq:3.10}
\end{eqnarray}%
where 
\begin{eqnarray}
F(\nabla \bm{u},\nabla \bm{\omega} ,\bm{\Omega}) &=&2\bm{\Omega R}(%
W_{\bm{E}})\cdot \nabla \bm{u}+\bm{\Omega R}(W_{\bm{\Gamma }%
})\cdot \nabla \bm{\omega }-\bm{\Omega R}(W_{\bm{E}})\bm{F}%
^{t}\cdot \bm{\Omega }+\bm{\bm{R}}^{t}\bm{\bm{\Omega F}}%
\cdot W_{\bm{EE}}[\bm{R}^{t}\bm{\Omega F}]  \nonumber \\
&&-2\bm{R}^{t}\nabla \bm{u}\cdot W_{\bm{EE}}[\bm{R}^{t}%
\bm{\Omega F}]-\bm{R}^{t}\nabla \bm{\omega }\cdot W_{\bm{%
\Gamma E}}[\bm{R}^{t}\bm{\Omega F}]-\bm{R}^{t}\bm{\Omega F}%
\cdot W_{\bm{E\Gamma }}[\bm{R}^{t}\nabla \bm{\omega}]. 
\label{eq:3.11}
\end{eqnarray}

If the equilibrium state is an energy minimizer, it is necessary that 
\begin{equation}
\ddot{E}\geq 0  \label{eq:3.12}
\end{equation}%
for all $\bm{u}$ and $\bm{\omega }$ such that $\bm{u}$ vanishes
on $\partial \kappa \setminus \partial \kappa _{t}$ and $\bm{\omega }$
vanishes on $\partial \kappa \setminus \partial \kappa _{c}.$


\section{The Legendre-Hadamard inequality}


\textit{Theorem}: If (\ref{eq:3.12}) is satisfied then it is necessary that the
Legendre-Hadamard inequality%
\begin{equation}
\bm{a}\otimes \bm{n}\cdot W_{\bm{EE}}[\bm{a}\otimes \bm{n}]+\bm{%
a}\otimes \bm{n}\cdot W_{\bm{E\Gamma }}[\bm{b}\otimes \bm{n}]+\bm{b}\otimes
\bm{n}\cdot W_{\bm{\Gamma E}}[\bm{a}\otimes \bm{n}]+\bm{b}\otimes \bm{n}\cdot
W_{\bm{\Gamma \Gamma }}[\bm{b}\otimes \bm{n}]\geq 0  \label{eq:4.1}
\end{equation}%
be satisfied at every $\bm{X}\in \kappa $ and for all vectors $\bm{%
a},\bm{b}$ and $\bm{n}.$

\qquad

\textit{Remark}: Choosing $\bm{a}$ or $\bm{b}$ to vanish in this
inequality yields the further necessary conditions%
\begin{equation}
\bm{a}\otimes \bm{n}\cdot W_{\bm{EE}}[\bm{a}\otimes \bm{n}]\geq 0\quad 
\textrm{and\quad }\bm{b}\otimes \bm{n}\cdot W_{\bm{\Gamma \Gamma }}[\bm{%
b}\otimes \bm{n}]\geq 0,  \label{eq:4.2}
\end{equation}%
again for every $\bm{X}\in \kappa $ and arbitrary $\bm{a},\bm{b}$ and $%
\bm{n}.$ Clearly these are also sufficient for (\ref{eq:4.1}) in the case of a
decoupled energy with $W_{\bm{E\Gamma }}=\bm{0}$ and $W_{\bm{%
\Gamma E}}=\bm{0}$. Further, (\ref{eq:4.1}) follows if $W$ is convex in the
strain measures $\bm{E}$ and $\bm{\Gamma }$ jointly. Indeed this
hypothesis underpins existence theorems for equilibria proved in \cite{neff2006existence,neff2015existence,lankeit2017integrability,neff2004geometrically} and guarantees that (\ref{eq:4.1}) is automatically satisfied at any
equilibrium state. However, this does not imply convexity of the overall
minimization problem due to the nonlinear nature of the strain measures. 

Moreover, in \cite{eremeyev2005acceleration,eremeyev2007constitutive} it is established that the strict form of the
Legendre-Hadamard inequality (\ref{eq:4.1}) ensures the propagation of acceleration
waves in dynamical Cosserat elasticity. Thus the classical connection
between this inequality and the reality of propagation speeds - well known
in the setting of conventional hyperelasticity - carries over to the
Cosserat framework.

\qquad

\textit{Proof of the Theorem}: Following \cite{giaquinta2004weierstrass} we consider variations%
\begin{equation}
\bm{u}(\bm{X})=\epsilon \bm{\xi} (\bm{Y})\textrm{\quad and\quad }\bm{\omega} (\bm{X})=\epsilon \bm{\eta} (\bm{Y})\textrm{\quad with\quad }\bm{Y}=\epsilon
^{-1}(\bm{X}-\bm{X}_{0}),  \label{eq:4.3}
\end{equation}%
where $\bm{X}_{0}$ is an interior point of $\kappa $, $\epsilon $ is a
positive constant, and $\bm{\xi} ,\bm{\eta}$ are compactly supported in a
region $D,$ the image of a strictly interior neighborhood $\kappa ^{\prime}\subset \kappa $ of $\bm{X}_{0}$ under the map $\bm{Y}(\bm{\cdot} ).$
Accordingly $\bm{u}$ and $\bm{\omega }$ (hence $\bm{\Omega }$)
vanish on $\partial \kappa $ and are therefore admissible. For these
variations (\ref{eq:3.12}) reduces, after dividing by $\epsilon ^{3}$, passing to the
limit $\epsilon \rightarrow 0$ and invoking the Dominated Convergence
Theorem, to%
\begin{equation}
\int_{D}\{\bm{R}_{0}^{t}\nabla \bm{\xi} \cdot \bm{\mathcal{A}}[\bm{R}%
_{0}^{t}\nabla \bm{\xi}]+\bm{R}_{0}^{t}\nabla \bm{\xi} \cdot 
\bm{\mathcal{B}}[\bm{R}_{0}^{t}\nabla \bm{\eta}]+\bm{R}_{0}^{t}\nabla 
\bm{\eta} \cdot \bm{\mathcal{B}}^{t}[\bm{R}_{0}^{t}\nabla \bm{\xi}]+%
\bm{R}_{0}^{t}\nabla \bm{\eta} \cdot \bm{\mathcal{C}}[\bm{R}_{0}^{t}\nabla
\bm{\eta}]\}dv\geq 0,  \label{eq:4.4}
\end{equation}%
where $\bm{R}_{0}=\bm{R}(\bm{X}_{0}),$ and with $\bm{\mathcal{A}}=W_{\bm{%
EE\mid X}_{0}}=\bm{\mathcal{A}}^{t},$ $\bm{\mathcal{B}}=W_{\bm{E\Gamma
\mid X}_{0}},$ $\bm{\mathcal{B}}^{t}=W_{\bm{\Gamma E\mid X}_{0}}$
and $\bm{\mathcal{C}}=W_{\bm{\Gamma \Gamma \mid X}_{0}}=\bm{%
\mathcal{C}}^{t}.$ Here and henceforth $\nabla $ is the gradient with
respect to $\bm{Y}$ and we have used the fact that $F,$ defined by
(\ref{eq:3.11}), vanishes in the limit.

We extend $\bm{\xi }$ and $\bm{\eta }$ to complex-valued vector
fields as 
\begin{equation}
\bm{\xi} =\bm{\xi}_{1}+i\bm{\xi}_{2}\textrm{\quad and\quad }\bm{\eta}
=\bm{\eta}_{1}+i\bm{\eta}_{2},  \label{eq:4.5}
\end{equation}%
where $\bm{\xi }_{1,2}$ and $\bm{\eta }_{1,2}$ are real-valued, and
use these to derive%
\begin{eqnarray}
&&\bm{R}_{0}^{t}\nabla \bm{\xi} \cdot \bm{\mathcal{B}}[\bm{R}_{0}^{t}\nabla 
\bm{\bar{\eta}}]+\bm{R}_{0}^{t}\nabla \bm{\eta} \cdot \bm{\mathcal{B}}%
^{t}[\bm{R}_{0}^{t}\nabla \bm{\bar{\xi}}]  \nonumber \\
&=&\bm{R}_{0}^{t}\nabla \bm{\xi}_{1}\cdot \bm{\mathcal{B}}[\bm{R}%
_{0}^{t}\nabla \bm{\eta }_{1}]+\bm{R}_{0}^{t}\nabla \bm{%
\xi }_{2}\cdot \bm{\mathcal{B}}[\bm{R}_{0}^{t}\nabla \bm{\eta}_{2}%
]+\bm{R}_{0}^{t}\nabla \bm{\eta}_{1}\cdot \bm{\mathcal{%
B}}^{t}[\bm{R}_{0}^{t}\nabla \bm{\xi}_{1}]+\bm{R}%
_{0}^{t}\nabla \bm{\eta}_{2}\cdot \bm{\mathcal{B}}^{t}[\bm{R}%
_{0}^{t}\nabla \bm{\xi }_{2}],  \label{eq:4.6}
\end{eqnarray}%
in which an overbar is used to denote the complex conjugate. The imaginary
part of this expression vanishes by virtue of the fact that $\bm{A}\cdot \bm{\mathcal{B}}[\bm{B}]=\bm{B}\cdot \bm{\mathcal{B}}^{t}[\bm{A}]$ for
arbitrary $\bm{A},\bm{B}.$ In the same way, we obtain%
\begin{eqnarray}
&&\bm{R}_{0}^{t}\nabla \bm{\xi} \cdot \bm{\mathcal{A}}[\bm{R}%
_{0}^{t}\nabla \bm{\bar{\xi}}]+\bm{R}_{0}^{t}\nabla \bm{\eta}
\cdot \bm{\mathcal{C}}[\bm{R}_{0}^{t}\nabla \bm{\bar{\eta}}]  \nonumber \\
&=&\bm{R}_{0}^{t}\nabla \bm{\xi }_{1}\cdot \bm{\mathcal{A}}%
[\bm{R}_{0}^{t}\nabla \bm{\xi }_{1}]+\bm{R}_{0}^{t}\nabla 
\bm{\xi }_{2}\cdot \bm{\mathcal{A}}[\bm{R}_{0}^{t}\nabla \bm{%
\xi }_{2}]+\bm{R}_{0}^{t}\nabla \bm{\eta }_{1}\cdot 
\bm{\mathcal{C}}[\bm{R}_{0}^{t}\nabla \bm{\eta }_{1}]+\bm{R}%
_{0}^{t}\nabla \bm{\eta }_{2} \cdot \bm{\mathcal{C}}[\bm{R}_{0}^{t}\nabla 
\bm{\eta }_{2}],  \label{eq:4.7}
\end{eqnarray}%
so that if (\ref{eq:4.3}) holds for real-valued $\bm{\xi }$ and $\bm{\eta},$
then it follows that%
\begin{equation}
\int_{D}\{\bm{R}_{0}^{t}\nabla \bm{\xi} \cdot \bm{\mathcal{A}}[\bm{R}%
_{0}^{t}\nabla \bm{\bar{\xi}}]+\bm{R}_{0}^{t}\nabla \bm{\xi}
\cdot \bm{\mathcal{B}}[\bm{R}_{0}^{t}\nabla \bm{\bar{\eta}}]+\bm{R}%
_{0}^{t}\nabla \bm{\eta} \cdot \bm{\mathcal{B}}^{t}[\bm{R}_{0}^{t}\nabla 
\bm{\bar{\xi}}]+\bm{R}_{0}^{t}\nabla \bm{\eta} \cdot \bm{\mathcal{C}}[\bm{R%
}_{0}^{t}\nabla \bm{\bar{\eta}}]\}dv\geq 0  \label{eq:4.8}
\end{equation}%
for complex-valued $\bm{\xi }$ and $\bm{\eta}.$

Consider%
\begin{equation}
\bm{\xi} (\bm{Y})=\bm{\alpha }\exp (ik\bm{n}\cdot\bm{Y})f(\bm{Y})%
\textrm{\quad and\quad }\bm{\eta} (\bm{Y})=\bm{\beta }\exp(ik\bm{%
n}\cdot \bm{Y})f(\bm{Y}),  \label{eq:4.9}
\end{equation}%
where $\bm{\alpha} ,\bm{\beta}$ and $\bm{n}$ are real fixed vectors, $k$
is a non-zero real number and $f$ is a real-valued differentiable function
compactly supported in $D.$ These yield%
\begin{equation}
\bm{R}_{0}^{t}\nabla \bm{\xi }=\exp (ik\bm{n}\cdot \bm{Y})(ikf\bm{%
a}\otimes \bm{n}+\bm{a}\otimes \nabla f)\quad \textrm{and\quad }\bm{R}%
_{0}^{t}\nabla \bm{\eta}=\exp(ik\bm{n}\cdot \bm{Y})(ikf\bm{b}\otimes
\bm{n}+\bm{b}\otimes \nabla f),  \label{eq:4.10}
\end{equation}%
with $\bm{a}=\bm{R}_{0}^{t}\bm{\alpha }$ and $\bm{b}=\bm{R}_{0}^{t}\bm{%
\beta}.$ Substitution into (\ref{eq:4.8})\ and division by $k^{2}$ results in%
\begin{eqnarray}
0 &\leq &\{\bm{a}\otimes \bm{n}\cdot \bm{\mathcal{A}}[\bm{a}\otimes \bm{n}]+2%
\bm{a}\otimes \bm{n}\cdot \bm{\mathcal{B}}[\bm{b}\otimes \bm{n}]+\bm{%
b}\otimes \bm{n}\cdot \bm{\mathcal{C}}[\bm{b}\otimes \bm{n}]\}\int_{D}f^{2}dv 
\nonumber \\
&&+k^{-2}\int_{D}\{\bm{a}\otimes \nabla f\cdot \bm{\mathcal{A}}[\bm{a}%
\otimes \nabla f]+2\bm{a}\otimes \nabla f\cdot \bm{\mathcal{B}}[%
\bm{b}\otimes \nabla f]+\bm{b}\otimes \nabla f\cdot \bm{\mathcal{%
C}}[\bm{b}\otimes \nabla f]\}dv.  \label{eq:4.11}
\end{eqnarray}%
Finally, as $k\rightarrow \infty $ we recover%
\begin{equation}
\bm{a}\otimes \bm{n}\cdot \bm{\mathcal{A}}[\bm{a}\otimes \bm{n}]+2\bm{a}\otimes \bm{n%
}\cdot \bm{\mathcal{B}}[\bm{b}\otimes \bm{n}]+\bm{b}\otimes \bm{n}\cdot 
\bm{\mathcal{C}}[\bm{b}\otimes \bm{n}]\geq 0,  \label{eq:4.12}
\end{equation}%
which is just (\ref{eq:4.1}) on account of the arbitrariness of $\bm{X}_{0}.$


\section{Example}


By way of illustration we apply inequalities (\ref{eq:4.2}) to the quadratic,
decoupled energy 
\begin{eqnarray}
W &=&\mu \left\Vert sym(\bm{E-I})\right\Vert ^{2}+\mu _{c}\left\Vert
skew(\bm{E-I})\right\Vert ^{2}+\frac{1}{2}\lambda \lbrack tr(\bm{%
E-I})\rbrack^{2}  \nonumber \\
&&+a_{1}\left\Vert dev(sym\bm{\Gamma})\right\Vert ^{2}+a_{2}\left\Vert
skew\bm{\Gamma}\right\Vert ^{2}+\frac{1}{3}a_{3}(tr\bm{\Gamma})%
^{2},  \label{eq:5.1}
\end{eqnarray}%
proposed in \cite{neff2017real} to model isotropic materials, where $\mu ,\mu
_{c},\lambda $ and $a_{1-3}$ are material constants and $\bm{I}$ is the
identity. Using the variational formulas $(\left\Vert \bm{A}\right\Vert
^{2})^{\cdot }=2\bm{A}\cdot \bm{\dot{A}}$ and $[(tr\bm{A})^{2}]^{\cdot
}=2(tr\bm{A})\bm{I}\cdot \bm{\dot{A}},$ together with the orthogonality of
symmetric and skew tensors, and also that of deviatoric and spherical
tensors, we obtain 
\begin{eqnarray}
\dot{W} &=&[2\mu sym(\bm{E-I})+2\mu _{c}skew(\bm{E-I})+\lambda tr(%
\bm{E-I})\bm{I}\rbrack\cdot \bm{\dot{E}}  \nonumber \\
&&+[2a_{1}dev(sym\bm{\Gamma})+2a_{2}skew\bm{\Gamma }+\frac{2}{3}%
a_{3}(tr\bm{\Gamma})\bm{I}\rbrack\cdot \bm{\dot{\Gamma}},  \label{eq:5.2}
\end{eqnarray}%
from which it follows that%
\begin{equation}
W_{\bm{E}}=2\mu sym(\bm{E-I})+2\mu _{c}skew(\bm{E-I})+\lambda tr(%
\bm{E-I})\bm{I}  \label{eq:5.3}
\end{equation}%
and 
\begin{equation}
W_{\bm{\Gamma }}=2a_{1}dev(sym\bm{\Gamma})+2a_{2}skew\bm{\Gamma 
}+\frac{2}{3}a_{3}(tr\bm{\Gamma})\bm{I}.  \label{eq:5.4}
\end{equation}

A further variation yields%
\begin{equation}
W_{\bm{EE}}[\bm{\dot{E}}]=2\mu sym\bm{\dot{E}}+2\mu _{c}skew%
\bm{\dot{E}}+\lambda (tr\bm{\dot{E}})\bm{I}  \label{eq:5.5}
\end{equation}%
and 
\begin{equation}
W_{\bm{\Gamma \Gamma }}[\bm{\dot{\Gamma}}]=2a_{1}dev(sym\bm{\dot{%
\Gamma}})+2a_{2}skew\bm{\dot{\Gamma}}+\frac{2}{3}a_{3}(tr\bm{\dot{%
\Gamma}})\bm{I}.  \label{eq:5.6}
\end{equation}%
Accordingly,%
\begin{equation}
W_{\bm{EE}}[\bm{a}\otimes \bm{n}]=\mu (\bm{a}\otimes \bm{n}+\bm{n}\otimes \bm{a})%
+\mu _{c}(\bm{a}\otimes \bm{n}-\bm{n}\otimes \bm{a})+\lambda (\bm{a}\cdot \bm{n})\bm{I} 
\label{eq:5.7}
\end{equation}%
and%
\begin{equation}
W_{\bm{\Gamma \Gamma }}[\bm{b}\otimes \bm{n}]=a_{1}[\bm{b}\otimes
\bm{n}+\bm{n}\otimes \bm{b}-\frac{2}{3}(\bm{b}\cdot \bm{n})\bm{I}]+a_{2}(\bm{b}\otimes
\bm{n}-\bm{n}\otimes \bm{b})+\frac{2}{3}a_{3}(\bm{b}\cdot \bm{n})\bm{I}.  \label{eq:5.8}
\end{equation}%
These in turn yield 
\begin{equation}
\bm{a}\otimes \bm{n}\cdot W_{\bm{EE}}[\bm{a}\otimes \bm{n}]=\mu \lbrack
\left\Vert \bm{a}\right\Vert ^{2}\left\Vert \bm{n}\right\Vert ^{2}+(%
\bm{a}\cdot \bm{n})^{2}]+\mu _{c}[\left\Vert \bm{a}\right\Vert
^{2}\left\Vert \bm{n}\right\Vert ^{2}-(\bm{a}\cdot \bm{n})^{2}]+\lambda (%
\bm{a}\cdot \bm{n})^{2}  \label{eq:5.9}
\end{equation}%
and%
\begin{equation}
\bm{b}\otimes \bm{n}\cdot W_{\bm{\Gamma \Gamma }}[\bm{b}\otimes \bm{n}]%
=a_{1}[\left\Vert \bm{b}\right\Vert ^{2}\left\Vert \bm{n}\right\Vert
^{2}+(\bm{b}\cdot \bm{n})^{2}]+a_{2}[\left\Vert \bm{b}\right\Vert
^{2}\left\Vert \bm{n}\right\Vert ^{2}-(\bm{b}\cdot \bm{n})^{2}]+\frac{2}{%
3}(a_{3}-a_{1})(\bm{b}\cdot \bm{n})^{2}.  \label{eq:5.10}
\end{equation}

Introducing angles $\alpha $ and $\beta $ defined by $\bm{a}\cdot \bm{n}%
=\left\Vert \bm{a}\right\Vert \left\Vert \bm{n}\right\Vert \cos
\alpha $ and $\bm{b}\cdot \bm{n}=\left\Vert \bm{b}\right\Vert \left\Vert 
\bm{n}\right\Vert \cos \beta $ we find that inequalities (\ref{eq:4.2}) are
satisfied if and only if%
\begin{eqnarray}
0 &\leq &\mu +\mu _{c}+(\mu -\mu _{c}+\lambda )\cos ^{2}\alpha  \nonumber \\
&=&(\mu +\mu _{c})(\cos ^{2}\alpha +\sin ^{2}\alpha )+(\mu -\mu _{c}+\lambda
)\cos ^{2}\alpha  \nonumber \\
&=&(2\mu +\lambda )\cos ^{2}\alpha +(\mu +\mu _{c})\sin ^{2}\alpha \textrm{%
\quad }  \label{eq:5.11}
\end{eqnarray}%
and%
\begin{eqnarray}
0 &\leq &a_{1}+a_{2}+[a_{1}-a_{2}+\frac{2}{3}(a_{3}-a_{1})]\cos ^{2}\beta 
\nonumber \\
&=&(a_{1}+a_{2})(\cos ^{2}\beta +\sin ^{2}\beta )+[a_{1}-a_{2}+\frac{2}{3}%
(a_{3}-a_{1})]\cos ^{2}\beta  \nonumber \\
&=&\frac{2}{3}(2a_{1}+a_{3})\cos ^{2}\beta +(a_{1}+a_{2})\sin ^{2}\beta , 
\label{eq:5.12}
\end{eqnarray}%
for all $\alpha $ and $\beta .$ The necessary and sufficient conditions 
\begin{equation}
2\mu +\lambda \geq 0,\quad \mu +\mu _{c}\geq 0,\quad 2a_{1}+a_{3}\geq 0\quad 
\textrm{and\quad }a_{1}+a_{2}\geq 0  \label{eq:5.12}
\end{equation}%
follow immediately, and coincide with the Legendre-Hadamard conditions
derived in \cite{neff2017real} for linearized, isotropic Cosserat elasticity.

We observe that in general (\ref{eq:4.1}) and (\ref{eq:4.2}) do not impose restrictions on the
constitutive function $W$, but rather on the configuration fields $\{\bm{%
\chi} (\bm{X}),\bm{R}(\bm{X})\}.$ In the present example, however, these emerge as
constitutive inequalities due to the quadratic nature of the energy (\ref{eq:5.1}).


\appendix
\section{Appendix}



To confirm the kinematic admissibility of the first and second variations $%
\bm{\dot{R}}$ and $\bm{\ddot{R}}$ defined by (\ref{eq:2.13})$_{3}$ and (\ref{eq:3.7}),
consider a tensor-valued function $\bm{Q}(\bm{X};\epsilon )$
satisfying the differential equation%
\begin{equation}
\bm{Q}^{\prime}=\bm{WQ}\textrm{\quad with\quad }\bm{Q}(\bm{X};%
0)=\bm{R}(\bm{X}),  \label{eq:A.1}
\end{equation}%
where $(\cdot )^{\prime }=\partial (\cdot )/\partial \epsilon ,$ $\bm{R}$
is a rotation, and $\bm{W}(\bm{X};\epsilon )$ is an arbitrary
differentiable skew tensor function. Let $\bm{Z}(\bm{X};\epsilon )=%
\bm{QQ}^{t}.$ Then,%
\begin{equation}
\bm{Z}^{\prime }=\bm{WZ}-\bm{ZW}\textrm{\quad with\quad }\bm{Z}(\bm{%
X};0)=\bm{I}.  \label{eq:A.2}
\end{equation}%
This has the unique solution $\bm{Z}(\bm{X};\epsilon )=\bm{I},$
implying that $\bm{Q}(\bm{X};\epsilon)$ is orthogonal with $\det 
\bm{Q}=\pm 1.$ Further,%
\begin{equation}
(\det \bm{Q})^{\prime }/\det \bm{Q}=tr(\bm{Q}^{\prime }\bm{Q}%
^{-1})=tr\bm{W}=0,  \label{eq:A.3}
\end{equation}%
implying that $\det \bm{Q}(\bm{X};\epsilon)=\det \bm{R}=1$ and
hence that $\bm{Q}(\bm{X};\epsilon )$ is an admissible Cosserat
rotation field. The notation $\bm{\dot{R}}=\bm{Q}_{\mid \epsilon =0}^{\prime
} $ then yields (\ref{eq:2.13})$_{3}$ with $\bm{\Omega} (\bm{X})=\bm{W}(\bm{X};0).$

From (\ref{eq:A.1}) we have 
\begin{equation}
\bm{Q}^{\prime \prime }=(\bm{WQ})^{\prime }=\bm{W}^{\prime }%
\bm{Q}+\bm{W}^{2}\bm{Q}  \label{eq:A.4}
\end{equation}%
in which $\bm{W}^{\prime }$ is skew. This integrates to 
\begin{equation}
\bm{Q}^{\prime }=\bm{WQ}+\bm{C}  \label{eq:A.5}
\end{equation}%
in which $\bm{C}$ independent of $\epsilon .$ Evaluating at $\epsilon =0$
yields $\bm{C}=\bm{\dot{R}}-\bm{\Omega R}$, which vanishes by (\ref{eq:2.13})$%
_{3}. $ Accordingly $\bm{Q}^{\prime }=\bm{WQ}$, which, as we have
seen, ensures that $\bm{Q}(\bm{X};\epsilon)$ is a rotation provided
that $\bm{R}=\bm{Q}(\bm{X};0)$ is a rotation. On setting $\epsilon =0$ in
(\ref{eq:A.4}) we recover (\ref{eq:3.7}) in which $\bm{\ddot{R}}=\bm{Q}_{\mid \epsilon
=0}^{\prime \prime }$ and $\bm{\Phi} = \bm{W}_{\mid \epsilon =0}^{\prime }$.
The arbitrariness of the skew function $\bm{W}(\bm{X};\epsilon )$ implies
that the skew tensor fields $\bm{\Omega} (\bm{X})$ and $\bm{\Phi} (\bm{X})$ in
(\ref{eq:2.13})$_{3}$ and (\ref{eq:3.7}) can be chosen independently and arbitrarily.

\qquad


\textit{Acknowledgments}: The work of MS and DJS was supported by the US NSF
through grant CMMI-1931064, and that of PN by the German DFG through grant
NE 902/8-1.

\bibliographystyle{unsrt}
\bibliography{reference.bib}

\begin{thebibliography}{10}

\bibitem{cosserat1909theorie}
E.~Cosserat and F.~Cosserat.
\newblock Th{\'e}orie des corps d{\'e}formables.
\newblock {\em Herman, Paris}, (1909).

\bibitem{toupin1964theories}
R.~A. Toupin.
\newblock Theories of elasticity with couple-stress.
\newblock {\em Arch. Ration. Mech. Anal.}, 17:85--112, (1964).

\bibitem{reissner1975note}
E.~Reissner.
\newblock Note on the equations of finite-strain force and moment stress
  elasticity.
\newblock {\em Stud. Appl. Math.}, 54(1):1--8, (1975).

\bibitem{reissner1987further}
E.~Reissner.
\newblock A further note on finite-strain force and moment stress elasticity.
\newblock {\em Z. angew. Math.Phys.}, 38(5):665--673, (1987).

\bibitem{truesdell2004non}
C.~Truesdell and W.~Noll.
\newblock The non-linear field theories of mechanics, 3rd edn. (ed. S.S.
  Antman), Springer, Berlin, (2004).

\bibitem{neff2006existence}
P.~Neff.
\newblock Existence of minimizers for a finite-strain micromorphic elastic
  solid.
\newblock {\em Proc. Roy. Soc. Edinburgh A}, 136(5):997--1012, (2006).

\bibitem{pietraszkiewicz2009natural}
W.~Pietraszkiewicz and V.A. Eremeyev.
\newblock On natural strain measures of the non-linear micropolar continuum.
\newblock {\em Int. J. Solids Structures}, 46(3-4):774--787, (2009).

\bibitem{neff2015existence}
P.~Neff, M.~B{\^\i}rsan, and F.~Osterbrink.
\newblock Existence theorem for geometrically nonlinear cosserat micropolar
  model under uniform convexity requirements.
\newblock {\em J. Elasticity}, 121(1):119--141, (2015).

\bibitem{lankeit2017integrability}
J.~Lankeit, P.~Neff, and F.~Osterbrink.
\newblock Integrability conditions between the first and second cosserat
  deformation tensor in geometrically nonlinear micropolar models and existence
  of minimizers.
\newblock {\em Z. angew. Math. Phys.}, 68(1):11, (2017).

\bibitem{neff2017real}
P.~Neff, A.~Madeo, G.~Barbagallo, M.~V. d'Agostino, R.~Abreu, and I.-. Ghiba.
\newblock Real wave propagation in the isotropic-relaxed micromorphic model.
\newblock {\em Proc. Roy. Soc. A}, 473(2197):20160790, (2017).

\bibitem{neff2004geometrically}
P.~Neff.
\newblock {\em Geometrically exact Cosserat theory for bulk behaviour and thin
  structures. Modelling and mathematical analysis.}
\newblock PhD thesis, Habil., TU-Darmstadt, (2004).

\bibitem{giaquinta2004weierstrass}
M.~Giaquinta and S.~Hildebrandt.
\newblock Calculus of Variations I, Springer, Berlin, (2004).

\bibitem{eugster2019continuum}
S.~Eugster, F.~dell'Isola, and D.~Steigmann.
\newblock Continuum theory for mechanical metamaterials with a cubic lattice
  substructure.
\newblock {\em Math. Mech. Compl. Sys.}, 7(1):75--98, (2019).

\bibitem{eremeyev2005acceleration}
Victor Eremeyev.
\newblock Acceleration waves in micropolar elastic media.
\newblock {\em Doklady Phys.}, 50(4):204--206, 2005.

\bibitem{eremeyev2007constitutive}
H.~Altenbach, V.~Eremeyev, L.~Lebedev, and L.A. Rend{\'o}n.
\newblock Acceleration waves and ellipticity in thermoelastic micropolar media.
\newblock {\em Arch. Appl. Mech.}, 80(3):217--227, 2010.

\end{thebibliography}

\end{document}